\documentclass[11pt]{amsart}

\usepackage{epsfig, amsmath, amssymb, amsthm, amsfonts}

\newcommand{\df}{\ensuremath{\partial}}

\newcommand{\alg}{\ensuremath{\mathcal{A}}}

\newcommand{\rr}{\ensuremath{\mathbb{R}}}
\newcommand{\zz}{\ensuremath{\mathbb{Z}}}

\newcommand{\mylabel}[1]{\label{#1}
}

\theoremstyle{plain}
\newtheorem{thm}{Theorem}[section]
\newtheorem{cor}[thm]{Corollary}
\newtheorem{lem}[thm]{Lemma}

\theoremstyle{definition}
\newtheorem{defn}[thm]{Definition}
\newtheorem*{pty}{Property}

\theoremstyle{remark}
\newtheorem*{rem}{Remark}
\newtheorem*{ex}{Example}

\begin{document}

\title{Augmentations and Rulings of Legendrian Knots}

\author[J. Sabloff]{Joshua M. Sabloff} 
\address{Haverford College, Haverford, PA 19041} 
\email{jsabloff@haverford.edu}

\begin{abstract}
  A connection between holomorphic and generating family invariants of
  Legendrian knots is established; namely, that the existence of a
  ruling (or decomposition) of a Legendrian knot is equivalent to the
  existence of an augmentation of its contact homology.  This result
  was obtained independently and using different methods by Fuchs and
  Ishkhanov \cite{fuchs-ishk}.  Close examination of the proof yields
  an algorithm for constructing a ruling given an augmentation.
  Finally, a condition for the existence of an augmentation in terms
  of the rotation number is obtained.
\end{abstract}


\maketitle

\section{Introduction}
\mylabel{sec:intro}

A fundamental problem in Legendrian knot theory is to determine when
two knots are (or are not) Legendrian isotopic.\footnote{See
  Section~\ref{sec:background} for the basic definitions in Legendrian
  knot theory.}  Bennequin's 1983 paper \cite{bennequin} started the
enterprise by introducing the two ``classical'' invariants of
Legendrian knots, the Thurston-Bennequin invariant $tb$ and the
rotation number $r$.  Classification results based on these invariants
followed in the early 1990's: Legendrian unknots \cite{yasha-fraser},
torus knots \cite{etnyre-honda:knots}, and figure eight knots
\cite{etnyre-honda:knots} are completely classified by their
topological type and classical invariants.

Starting in the late 1990's, two methods for constructing
``non-classical'' invariants of Legendrian knots were developed. The
first is a relative version of the contact homology of Eliashberg,
Givental, and Hofer \cite{egh}. This theory uses holomorphic
techniques to associate a non-commutative differential graded algebra
(DGA) to a knot diagram. The homology of the DGA is invariant under
Legendrian isotopy.  In \cite{chv}, Chekanov rendered this theory
combinatorially computable.  He then used a linearized version of it
to distinguish examples of Legendrian $5_2$ knots in the standard
contact $\rr^3$ that have the same classical invariants.\footnote{This
  invariant is also referred to as the Chekanov-Eliashberg DGA in the
  literature.}  The second method is based on generating families,
i.e.\ families of functions whose critical values generate fronts of
Legendrian knots.  Chekanov has produced ``ruling'' or
``decomposition'' invariants based on generating families that can
distinguish his original $5_2$ examples (see \cite{chv:survey,
  chv-pushkar}).  In addition, Traynor has fashioned a non-classical
theory based on generating families for Legendrian links in the solid
torus \cite{lisa:links}.

The goal of this paper is to strengthen a connection, discovered by
Fuchs \cite{fuchs:augmentations}, between the ability to linearize the
contact homology DGA and the non-vanishing of Chekanov's count of
rulings for Legendrian knots in the standard contact $\rr^3$.  It is
possible to linearize the contact homology DGA if and only if there
exists an augmentation, i.e.\ a map $\varepsilon$ from the algebra to
the base ring that sends the image of the differential to zero.  It is
useful to further stipulate that the augmentation has support on
generators of grading zero modulo $2r(K)$ or grading divisible by a
divisor $\rho$ of $2r(K)$.  The former are ``graded'' augmentations;
the latter are ``$\rho$-graded'' augmentations.  Fuchs' original
result was:

\begin{thm}[Fuchs \cite{fuchs:augmentations}]
  \mylabel{thm:fuchs-augm} If a front diagram of a Legendrian knot $K$
  has a (graded or $\rho$-graded) normal ruling, then the contact
  homology DGA of $K$ has a (graded or $\rho$-graded) augmentation.
\end{thm}

The central result of this paper is the converse, which Fuchs and
Ishkhanov have proved independently, using different methods, in
\cite{fuchs-ishk}:

\begin{thm}
  \mylabel{thm:main} If the contact homology DGA of a Legendrian knot
  $K$ has a (graded or $\rho$-graded) augmentation, then any front
  diagram of $K$ has a (graded or $\rho$-graded) normal ruling.
\end{thm}

A consequence is an easy criterion for checking if the contact
homology DGA of a Legendrian knot has an augmentation:

\begin{thm}
  \mylabel{thm:rot-ruling} If the Chekanov-Eliashberg DGA of a
  Legendrian knot $K$ has a $2$-graded augmentation, then its rotation
  number is zero.
\end{thm}

These results contribute to recent work that examines the relationship
between the contact homology and generating family approaches to
constructing non-classical invariants. Ng and Traynor found that a
linearized version of the contact homology DGA and generating family
homology contain the same information for a large class of
two-component links in the solid torus \cite{lenny-lisa}. Work of Zhu
\cite{zhu} and of Ekholm, Etnyre, and Sullivan \cite{ees:graph-trees}
(see also \cite{fukaya-oh}) shows that a different sort of generating
function homology that uses ``graph trees'' for a single set of
generating functions can be used to compute the contact homology DGA.
The ideas behind this work have already provided the motivation for
Ng's combinatorial construction of invariants of topological braids
and knots using the contact homology of Legendrian tori in $ST^*\rr^3$
\cite{lenny:knot-invts-1}.

The rest of the paper is organized as follows:
Section~\ref{sec:background} lays out the necessary background and
notation for diagrams of Legendrian knots, the contact homology DGA,
and normal rulings.  Section~\ref{sec:main} contains the proof of
Theorem~\ref{thm:main} using a modification of a plat diagram of a
Legendrian knot.  Finally, the proof of Theorem~\ref{thm:rot-ruling}
appears in Section~\ref{sec:augm-r}.

\subsection*{Acknowledgments}

This paper has greatly benefited from discussions with John Etnyre,
Lisa Traynor, Lenny Ng, and especially Paul Melvin, who first
conjectured Theorem~\ref{thm:rot-ruling} based on computations done by
himself and Sumana Shrestha.  Conversations with Dmitry Fuchs helped
to clarify the hypotheses in Theorem~\ref{thm:rot-ruling}, and the
referee's comments greatly improved the exposition of the proof of
Theorem~\ref{thm:main}.

\section{Background Notions}
\mylabel{sec:background}

\subsection{Diagrams of Legendrian Knots}
\mylabel{ssec:knots}

This section briefly reviews some basic notions of Legendrian knot
theory.  For a more comprehensive introduction, see
\cite{etnyre:knot-intro,lecnotes}.

The \textbf{standard contact structure} on $\rr^3$ is the completely
non-integrable $2$-plane field given by the kernel of $\alpha = dz -
y\,dx$.   A
\textbf{Legendrian knot} is an embedding $K: S^1 \to \rr^3$ that is
everywhere tangent to the contact planes.  In particular, the
embedding must satisfy
\begin{equation} \mylabel{eqn:leg-cond}
  \alpha(K') = 0.
\end{equation}
An ambient isotopy of $K$ through Legendrian knots is a
\textbf{Legendrian isotopy}.

There are two useful projections of Legendrian knots.  The
\textbf{Lagrangian projection} is given by the map
\begin{equation*}
  \pi_l: (x,y,z) \mapsto (x,y),
\end{equation*}
while the \textbf{front projection} is given by
\begin{equation*}
  \pi_f: (x,y,z) \mapsto (x,z).
\end{equation*}
The Lagrangian and front projections of a Legendrian trefoil knot
appear in Figure~\ref{fig:numbered-trefoil}.  

\begin{figure}
  \centerline{\includegraphics{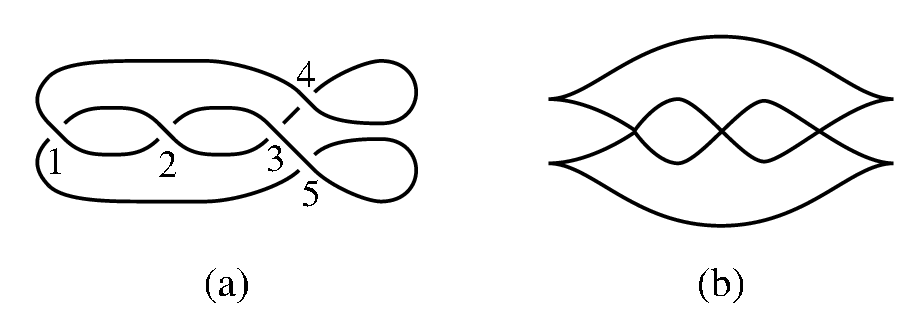}}
  \caption{(a) Lagrangian and (b) front diagrams of a Legendrian trefoil
    knot.  The meaning of the numbers in (a) will become clear in
    Section~\ref{ssec:chv-el}.}  \mylabel{fig:numbered-trefoil}
\end{figure}

In the front projection, the $y$ coordinate of a knot may be recovered
from the slope of the front projection via (\ref{eqn:leg-cond}):
\begin{equation} \mylabel{eqn:y-slope}
  y = \frac{dz}{dx}.
\end{equation}
This fact has several consequences:
\begin{itemize}
\item The front projection of a Legendrian knot is never vertical.
  Instead of vertical tangencies, front projections have \emph{cusps}
  like those on the extreme left and right of
  Figure~\ref{fig:numbered-trefoil}(b).
\item There is no need to specify crossing information at a double
  point: the strand with the smaller slope always has a smaller $y$
  coordinate.  This means that it will pass in \emph{front} of the
  strand with the larger slope, as the $y$ axis must point into the
  page in the front projection.
\item Any circle in the $xz$ plane that has no vertical tangencies and
  that is immersed except at finitely many cusps lifts to a Legendrian
  knot via equation (\ref{eqn:y-slope}).
\end{itemize}
A front diagram is in \textbf{plat position} if all of the left cusps
have the same $x$ coordinate, all of the right cusps have the same $x$
coordinate, and no two crossings have the same $x$ coordinate. For
example, the diagram of the trefoil in
Figure~\ref{fig:numbered-trefoil}(b) is in plat position.  The $x$
coordinates of the crossings and cusps are the \textbf{singular
  values} of the front.  Any front diagram may be put into plat
position using Legendrian versions of Reidemeister type II moves and
planar isotopy.

Though the front projection is easier to work with, it is more natural
to define the contact homology DGA using the Lagrangian projection.
Ng's \textbf{resolution} procedure (see \cite{lenny:computable}, and
Figures 2 and 3 in particular) gives a canonical translation from a
front diagram to a Lagrangian diagram.  This procedure, in fact, was
used to derive the Lagrangian projection in
Figure~\ref{fig:numbered-trefoil}(a) from the front projection in
Figure~\ref{fig:numbered-trefoil}(b).  Combinatorially, there are
three steps:
\begin{enumerate}
\item Smooth the left cusps;
\item Replace the right cusps with a loop (see the right side of the
  Lagrangian projection in Figure~\ref{fig:numbered-trefoil}); and
\item Resolve the crossings so that the overcrossing strand is the one
  with smaller slope.
\end{enumerate}
A key feature of the resolution procedure is that the heights of
the crossings in the Lagrangian projection strictly increase from left
to right, with the jumps in height between crossings as large as
desired.

As mentioned in the introduction, there are two ``classical''
invariants for Legendrian knots up to Legendrian isotopy.  The first
classical invariant is the \textbf{Thurston-Bennequin number} $tb(K)$,
which measures the twisting of the contact planes around the knot $K$.
The second classical invariant, the \textbf{rotation number} $r(K)$,
is defined for \emph{oriented} Legendrian knots.  It measures the
turning of the tangent direction to $K$ inside the contact planes with
respect to the trivialization given by the vector fields $\partial_y$
and $\partial_x + y\partial_z$.  The rotation number of an oriented
Legendrian knot $K$ may be computed using the rotation number of the
tangent vector to the Lagrangian projection in the plane.  In the
front projection, the rotation number is half of the difference
between the number of downward-pointing cusps and the number of
upward-pointing cusps.

\subsection{The Contact Homology DGA and Augmentations}
\mylabel{ssec:chv-el}

This section contains a brief review of the definition of the contact
homology DGA of a Legendrian knot.  The DGA was originally defined by
Chekanov in \cite{chv} for Lagrangian diagrams; see also \cite{ens}.

Let $K$ be an oriented Legendrian knot in the standard contact $\rr^3$
with a generic Lagrangian diagram $\pi_l(K)$.  Label the crossings by
$q_1, \ldots, q_n$.  Let $\alg$ be the graded free unital tensor
algebra over $\zz/2$ generated by the set $\{q_1, \ldots,
q_{n}\}$.\footnote{It is possible to define the algebra over $\zz[T,
  T^{-1}]$; see \cite{ens}.}  To define the grading, a capping path
$\gamma_i$ needs to be assigned to each crossing.  A \textbf{capping
  path} is one of the two paths in $\pi_l(K)$ that starts at the
overcrossing of $q_i$ and ends when $\pi_l(K)$ first returns to $q_i$,
necessarily at an undercrossing.  Assume, without loss of generality,
that the strands of $\pi_l(K)$ at each crossing are orthogonal. The
\textbf{grading} of $q_i$ is:
\begin{equation*}
  |q_i| \equiv 2 r(\gamma_i) - \frac{1}{2} \pmod{2r(K)}.
\end{equation*}
Extend the grading to all words in $\alg$ by letting the grading of a
word be the sum of the gradings of its constituent generators.

\begin{rem}
  It is simple to assign gradings directly from a plat diagram. Assign
  a grading of $1$ to each generator coming from a right cusp.  To
  assign a grading to a crossing, begin as in \cite{chv:survey} by
  letting $C(K)$ be the set of points on $K$ corresponding to cusps of
  $\pi_f(K)$.  The \textbf{Maslov index} is a locally constant
  function
  \begin{equation*}
    \mu: K \setminus C(K) \to \zz/2r(K)
  \end{equation*}
  that satisfies the relations depicted in Figure~\ref{fig:maslov}
  near the cusps.  This function is well-defined up to an overall
  constant.  Near a crossing $q_i$, let $\alpha_i$ (resp.  $\beta_i$)
  be the strand of $\pi_f(K)$ with more negative (resp.  positive)
  slope.  Assign the grading $|q_i| \equiv \mu(\alpha_i) - \mu(\beta_i)
  \pmod{2r(K)}$.
\end{rem}

\begin{figure}
  \centerline{\includegraphics{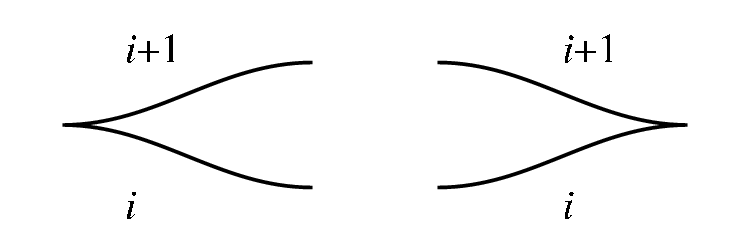}}
  \caption{The function $\mu$ must satisfy these relations near the 
    cusps.}  \mylabel{fig:maslov}
\end{figure}

The next step is to define a differential on $\alg$ by counting
certain immersions of the disk into $\pi_l(K)$.  Label the corners of
$\pi_l(K)$ as in Figure~\ref{fig:disks}(a).  The immersions of
interest are the following:

\begin{defn}
  \mylabel{defn:ce-immersed} Given a generator $q_i$ and an ordered
  set of generators $\mathbf{q} = \{ q_{j_1}, \ldots, q_{j_k}\}$, let
  $\Delta(q_i; \mathbf{q})$ be the set of orientation-preserving
  immersions
  \begin{equation*}
    f: D^2 \to \rr^2
  \end{equation*}
  that map $\partial D^2$ to $\pi_l(K)$ (up to smooth
  reparametrization), with the property that the restriction of $f$ to
  the boundary is an immersion except at the points $q_i, q_{j_1},
  \ldots, q_{j_k}$ and these points are encountered in
  counter-clockwise order along the boundary.  In a neighborhood of
  $q_i$ and the points in $\mathbf{q}$, the image of the disk under
  $f$ has the form indicated in Figure~\ref{fig:disks}(b) near $q_i$
  and in Figure~\ref{fig:disks}(c) near $q_{j_l}$.
\end{defn}

\begin{figure}
  \centerline{\includegraphics{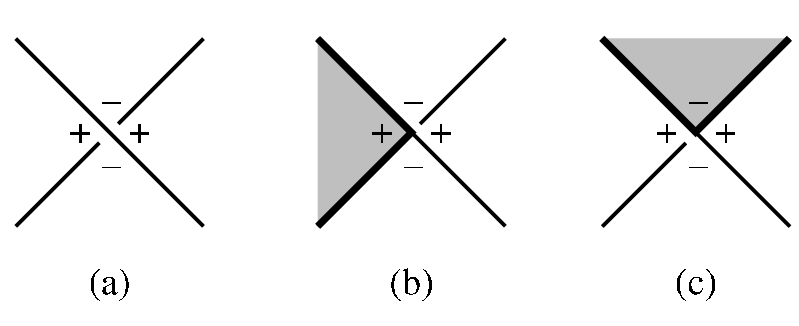}}
  \caption{(a) A labeling of the quadrants surrounding a crossing; (b)    
    the image of $f \in \Delta(q_i; \mathbf{q})$ near the crossing
    $q_i$ (the other $+$ quadrant is also possible); and (c) the image
    of $f$ near the crossing $q_{j_l}$ (the other $-$ quadrant is also
    possible).}  \mylabel{fig:disks}
\end{figure}

Finally, define the differential as follows:

\begin{defn}
  \mylabel{defn:ce-diffl} The differential \df\ is defined on a
  generator $q_i$ by the formula:
  \begin{equation}
    \df q_i = \sum_{\Delta(q_i; \mathbf{q})} \#\left(
    \Delta(q_i; \mathbf{q}) \right) q_{j_1} \cdots
    q_{j_k},
  \end{equation}
  where $\# \Delta(q_i; \mathbf{q})$ is the number of elements in the
  set $\Delta(q_i; \mathbf{q})$, counted modulo 2.  Extend $\df$ to
  all of $\alg$ via linearity and the Leibniz rule.
\end{defn}

Note that the sum in the definition of $\df$ is finite, and that if
$\Delta(q_i; \mathbf{q})$ is nonempty, then the height of the crossing
at $q_i$ is greater than the sum of the heights of the crossings
$q_{j_l}$; see \cite{chv}.

\begin{rem}
  In a diagram coming from the resolution of a plat diagram, the disks
  in the differential take on a simple form:
  \begin{enumerate}
  \item The disks are embedded, and 
  \item The intersection of any vertical line with a disk is
    connected.
  \end{enumerate}
\end{rem}

\begin{ex}
  Number the crossings of the trefoil knot as in
  Figure~\ref{fig:numbered-trefoil}.  The first three crossings have
  grading $0$, whereas the crossings that come from cusps in the plat
  diagram have grading $1$.  The only nontrivial differentials are:
  \begin{equation*}
    \begin{split}
      \df q_4 &= 1 + q_1 + q_3 + q_1 q_2 q_3, \\
      \df q_5 &= 1 + q_1 + q_3 + q_3 q_2 q_1. 
    \end{split}
  \end{equation*}
\end{ex}

The central results in this theory are:

\begin{thm}[\cite{chv}] 
  \mylabel{thm:dga}
  \begin{enumerate}
  \item The differential \df\ has degree $-1$.
  \item The differential satisfies $\df^2=0$.
  \item The ``stable tame isomorphism class'' of the DGA is invariant
    under Legendrian isotopy.
  \end{enumerate}
\end{thm}

The ``stable'' in part (3) of the theorem comes from the following
operation on a DGA $(\alg, \df)$: the \textbf{degree $i$
  stabilization} $S_i(\alg, \df)$ adds two new generators $\beta$ and
$\alpha$ to the algebra, where
\begin{equation*}
  |\beta| = i \text{ and } |\alpha| = i-1,
\end{equation*}
and the differential is extended to the new generators by:
\begin{equation*}
  \df \beta = \alpha \text{ and } \df \alpha = 0.
\end{equation*}
For the purposes of this paper, a stable tame isomorphism between two
DGAs $(\alg, \df)$ and $(\alg', \df')$ is a DGA isomorphism
\begin{equation*}
  \psi: S_{i_1} \left( \cdots S_{i_m}(\alg) \cdots \right) \to
  S_{j_1} \left( \cdots S_{j_n}(\alg') \cdots \right).
\end{equation*}
 
It is not easy to use the DGA to distinguish between Legendrian knots,
as it --- and its homology --- are fairly complicated objects.
Chekanov found computable invariants by linearizing the DGA.  Asking
whether the DGA has a graded augmentation is a first step in
generating linearized invariants:

\begin{defn}
  \mylabel{defn:augmentation} An \textbf{augmentation} is an algebra
  map $\varepsilon: \alg \to \zz/2$ that satisfies $\varepsilon \circ
  \df = 0$ and $\varepsilon(1) = 1$.  If, in addition, the
  augmentation has support on generators of degree zero, then it is
  \textbf{graded}; if it has support on generators divisible by a
  divisor $\rho$ of $2r(K)$, then it is $\rho$-\textbf{graded}.
\end{defn}

It is easy to extend a (graded or $\rho$-graded) augmentation over a
stabilization: simply send both $\beta$ and $\alpha$ to $0$.  In the
case of a degree $0$ stabilization --- or degree divisible by $\rho$
in the $\rho$-graded case --- there is another possible extension:
\begin{equation*}
  \varepsilon (\beta) = 1 \text{ and } \varepsilon(\alpha) = 0.
\end{equation*}
That is, if $|\beta|=0$, $\varepsilon(\beta)$ can be either $0$ or
$1$.  Either way, Theorem~\ref{thm:dga}(3) implies:

\begin{cor}
  \mylabel{cor:augm-invt} The existence of a (graded or $\rho$-graded)
  augmentation is invariant under Legendrian isotopy.
\end{cor}

\begin{ex}
  The DGA for the trefoil knot in the previous example has five graded
  augmentations.  For grading reasons, all of the augmentations are
  zero on $q_4$ and $q_5$, and it is easy to check that the following
  assignments work:
  \begin{center}
    \begin{tabular}{r||c|c|c}
      & $q_1$ & $q_2$ & $q_3$ \\ \hline
      $\varepsilon_1$ & 1 & 0 & 0 \\
      $\varepsilon_2$ & 1 & 1 & 0 \\
      $\varepsilon_3$ & 1 & 1 & 1 \\
      $\varepsilon_4$ & 0 & 1 & 1 \\
      $\varepsilon_5$ & 0 & 0 & 1
    \end{tabular}
  \end{center}
\end{ex}

\subsection{Rulings}
\mylabel{ssec:rulings}

The other object involved in Theorem~\ref{thm:main} is a (graded or
$\rho$-graded) normal ruling.  Suppose that a Legendrian knot $K$ has
a front diagram whose singular values all have distinct $x$
coordinates.  A \textbf{ruling} of such a front diagram of $K$
consists of a one-to-one correspondence between the set of left cusps
and the set of right cusps and, for each pair of corresponding cusps,
two paths in the front diagram that join them.  The ruling paths must
satisfy the following conditions:
\begin{enumerate}
\item Any two paths in the ruling meet only at crossings or at cusps;
  and
\item The interiors of the two paths joining corresponding cusps are
  disjoint, and hence they meet only at the cusps and bound a
  topological disk. Note that these disks are similar to those used to
  define the differential \df, but they may have ``obtuse'' corners;
  see Figure~\ref{fig:switch-config}(b), for example.
\end{enumerate}
As Fuchs notes, these conditions imply that the paths cover the front
diagram and the $x$ coordinate of each path in the ruling is
monotonic. 

At a crossing, either the two ruling paths incident to the crossing
pass through each other or one path lies entirely above the other.  In
the latter case, say that the ruling is \textbf{switched} at the
crossing.  Near a crossing, call the two ruling paths that intersect
the crossing \textbf{crossing paths} and the ruling paths that are
paired with the crossing paths \textbf{companion paths}.  If all of
the switched crossings of a ruling are of types (a--c) in
Figure~\ref{fig:switch-config}, then the ruling is \textbf{normal}.
If all of the switched crossings have grading $0$ (resp. grading
divisible by $\rho$), then the ruling is \textbf{graded} (resp.
$\rho$-graded).  It is not hard to see that in a graded ruling, both
crossing paths have the same Maslov index in configurations (a--c), as
do the companion paths in configurations (b) and (c).

\begin{figure}
  \centerline{\includegraphics{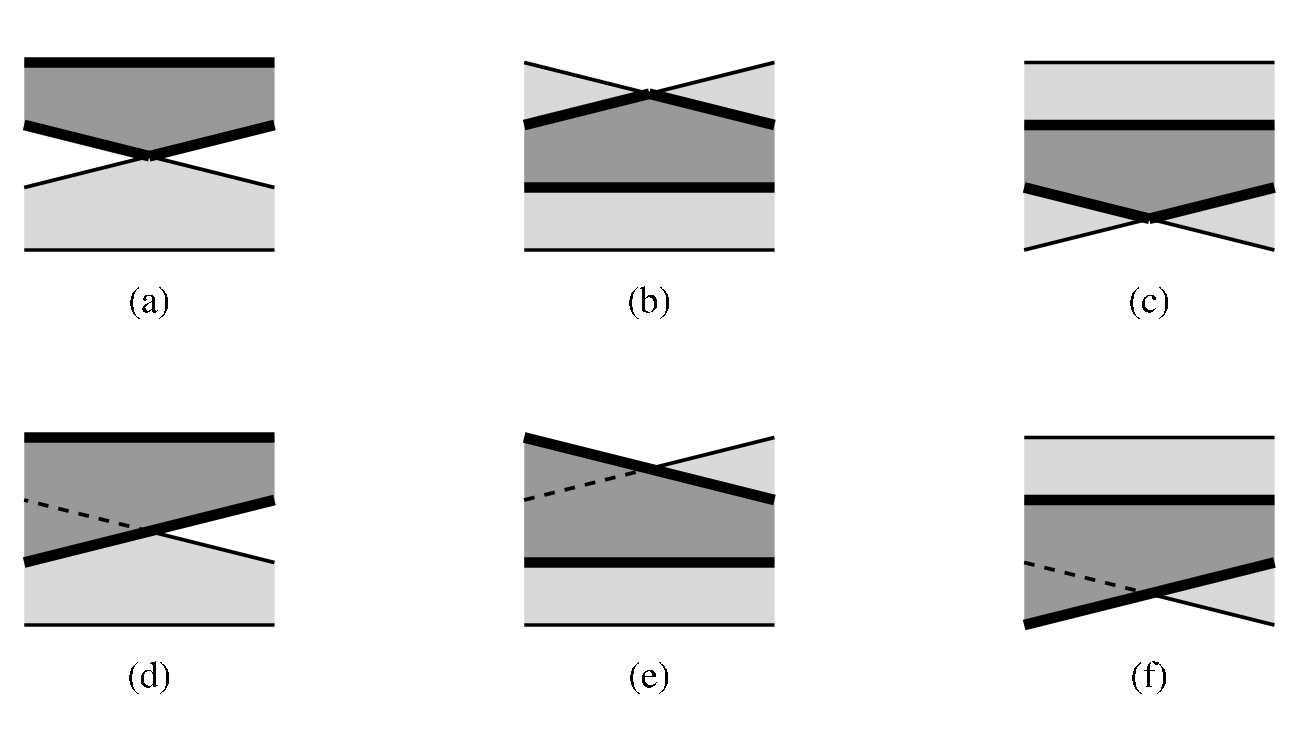}}
  \caption{The possible configurations of a normal ruling near a
    crossing.  Only the crossing paths and their companions are shown.
    Reflections of configurations (d--f) through a vertical axis are
    also allowed.}  \mylabel{fig:switch-config}
\end{figure}

\begin{ex}
  The trefoil pictured in Figure~\ref{fig:numbered-trefoil} has
  exactly three graded normal rulings.  They are pictured in
  Figure~\ref{fig:trefoil-rulings}.
\end{ex}

\begin{figure}
  \centerline{\includegraphics{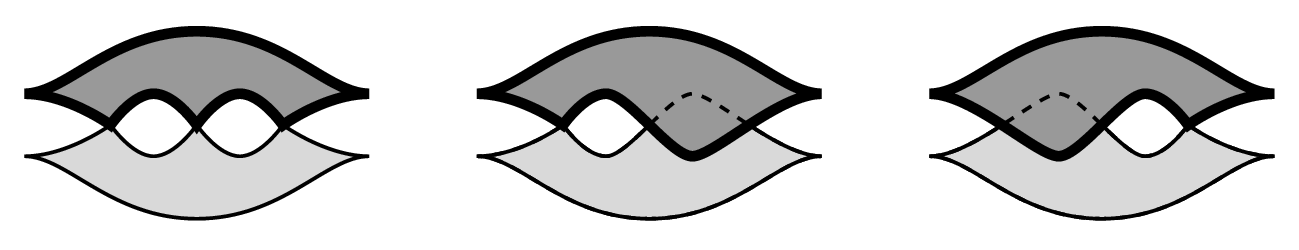}}
  \caption{The three normal rulings of the trefoil knot in 
    Figure~\ref{fig:numbered-trefoil}.}  \mylabel{fig:trefoil-rulings}
\end{figure}

The following theorem of Chekanov shows that normal rulings are
interesting objects in Legendrian knot theory:

\begin{thm}[Chekanov \cite{chv:survey}] \mylabel{thm:ruling-invt}
  The number of (graded or $\rho$-graded) normal
  rulings\footnote{Chekanov calls them \textbf{admissible
      decompositions} in \cite{chv:survey}; Chekanov and Pushkar call
    them \textbf{positive involutions} in \cite{chv-pushkar}.} is
  invariant under Legendrian isotopy.
\end{thm}

\section{From Augmentation to Ruling}
\mylabel{sec:main}

In light of Corollary~\ref{cor:augm-invt} and
Theorem~\ref{thm:ruling-invt}, the proof of Theorem~\ref{thm:main} ---
that the existence of an augmentation implies the existence of a
ruling --- only needs to consider Lagrangian diagrams that come from
resolving plats. The proof consists of extending the ruling crossing
by crossing from left to right.  The extension procedure will produce
only (graded or $\rho$-graded) normal switches, so the challenge will
be to prove that the paths paired in the ruling match up at the right
cusps.  To do this, the proof adopts Fuchs' philosophy of using
Legendrian isotopy to simplify the differential at the expense of
expanding the number of generators.  In practice, this means
converting plat diagrams into ``dipped diagrams'' in which certain
crossings are closely related to rulings (see
Section~\ref{ssec:dipped}).  By tracing the original augmentation
through the stable tame isomorphisms that relate the DGAs of the
original diagram and of the dipped diagram (see
Section~\ref{ssec:type2}), it will be possible to use properties of
the augmentation of the dipped diagram to conclude that the ruling
paths match at the right cusps (see Section~\ref{ssec:tracing}).

\subsection{Dipped Diagrams}
\mylabel{ssec:dipped}

A \textbf{dip} in a plat diagram looks innocent in the front
projection: it appears as the small wiggles pictured in
Figure~\ref{fig:dips}(a).  The new front is clearly isotopic to the
original one.  The Lagrangian diagram, however, has changed
dramatically; see Figure~\ref{fig:dips}(b).

\begin{figure}
  \centerline{\includegraphics{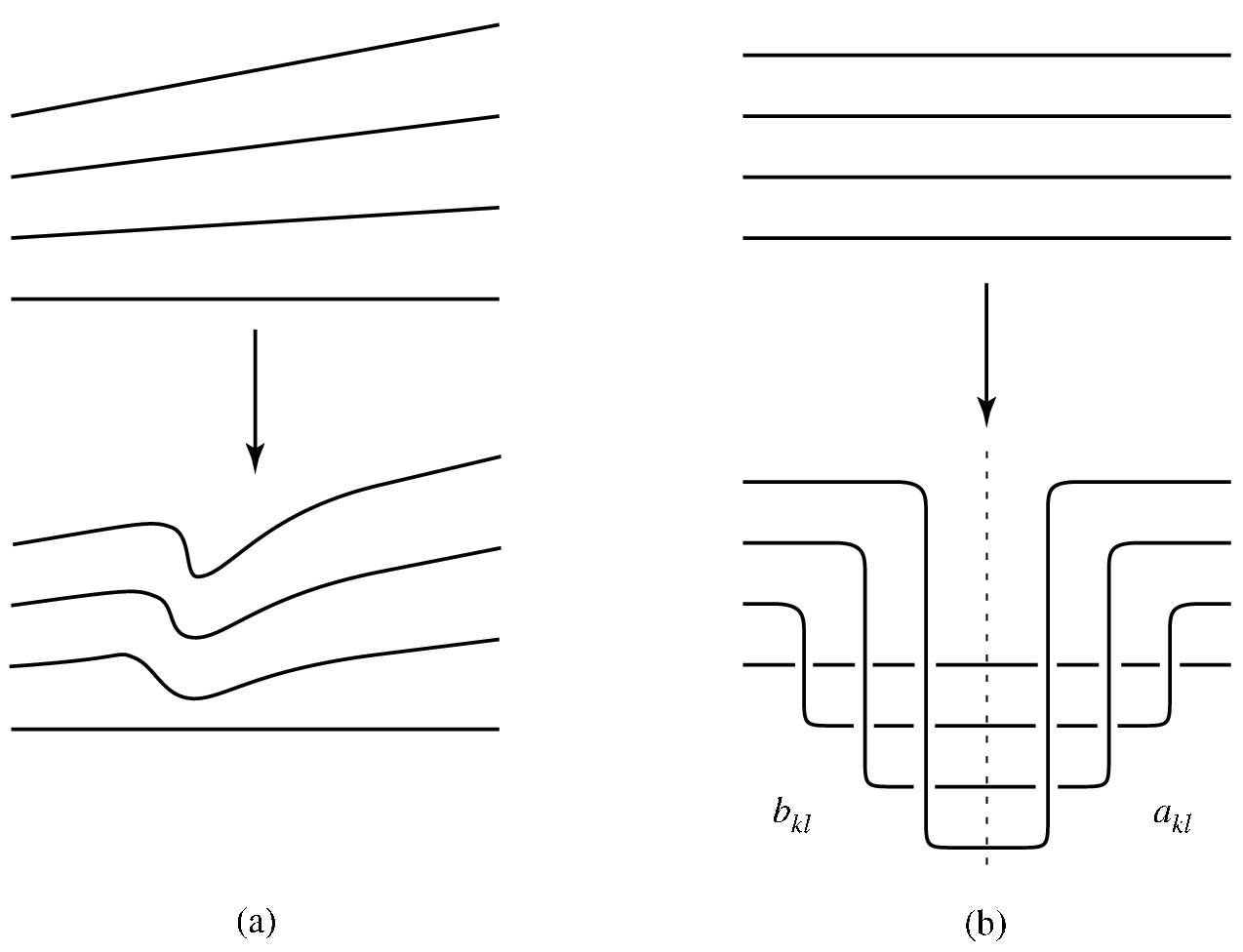}}
  \caption{(a) The modification in the front projection; (b) The
    modification in the Lagrangian projection (after a small planar
    isotopy).}  \mylabel{fig:dips}
\end{figure}

To see the transition to the dipped diagram in the Lagrangian
projection in terms of Reidemeister moves, start by numbering the
strands from bottom to top.  Using a Type II move, push strand $k$
over strand $l$ ($k>l$) in ascending lexicographic order, e.g.\ $3$
crosses $2$ after $3$ crosses $1$, and $4$ crosses $1$ after $3$
crosses $2$. If $k$ crosses $l$ after $i$ crosses $j$, write $(i,j)
\prec (k,l)$.  The new generators for the modified diagram are simple
to describe: assuming $k > l$, denote by $b_{kl}$ the leftmost
crossing of the strands $k$ and $l$ and by $a_{kl}$ the rightmost
crossing.  Say that the $b_{kl}$ generators belong to the
$b$-\textbf{lattice} and the $a_{kl}$ generators belong to the
$a$-\textbf{lattice}; see Figure~\ref{fig:dips}(b). It is not hard to
check that $|b_{kl}| = \mu(l) - \mu(k)$; note that this is the
negative of the grading of a crossing $q_i$.  Since the differential
lowers degree by $1$, it follows that $|a_{kl}| = |b_{kl}| - 1$.

The differential interacts straightforwardly with the new generators:

\begin{lem} \mylabel{lem:dip-disks}
  Suppose that $a$ and $b$ are the new crossings created by a Type II
  move during the creation of a dip, and let $y$ be any other
  crossing.  The generator $a$ appears at most once in any term of
  $\df y$, and if $a$ appears in $\df y$, then $b$ does not.
\end{lem}

\begin{proof}
  Consider a disk with a negative corner at $a$.  As shown in
  Figure~\ref{fig:partial-dip}, this corner must lie in the bottom
  left or top right quadrant adjacent to $a$.  In the case where the
  corner is at the bottom left, there is only one possible disk that
  comes from $\df b$.  Otherwise, the corner is at the top right and
  there are two cases.  First, suppose that the next corner on the
  upper strand lies in the $a$-lattice.  The disk must then lie
  entirely inside the $a$-lattice, as pictured in
  Figure~\ref{fig:partial-dip}(a). In particular, the disk satisfies
  the conditions of the lemma.
  
  Second, suppose that the next corner on the upper strand lies
  outside --- and hence to the right of --- the $a$-lattice. As shown
  in Figure~\ref{fig:partial-dip}(b), the lower strand must also exit
  the $a$-lattice without any further corners.  Note that the dipped
  diagram comes from modifying the resolution of a ``simple'' front
  (see \cite{lenny:computable}, Section 2.3) whose right cusps are all
  pushed out to the right, so any portion of a disk lying to the right
  of the dip must have connected vertical slices.  It follows that the
  rest of the disk must lie to the right of the figure, and hence that
  the disk satisfies the conditions of the lemma.
\end{proof}

  \begin{figure}
    \centerline{\includegraphics{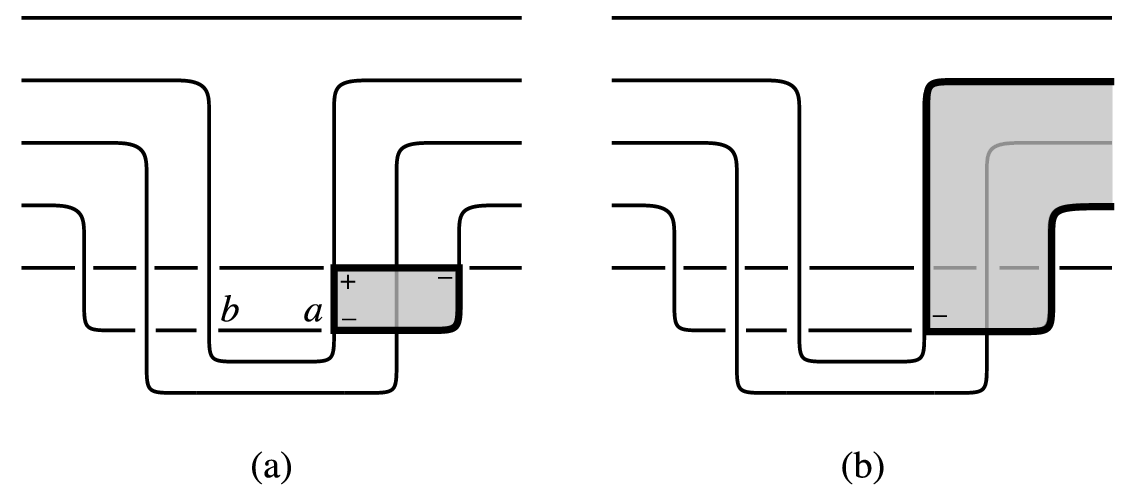}}
    \caption{The form of disks that have a negative corner at $a$
      after a Type II move.}
    \mylabel{fig:partial-dip}
  \end{figure}

\subsection{Type II Moves and DGA Maps}
\mylabel{ssec:type2}

In order to understand how the augmentations before and after the
formation of a dip are related, a closer examination of the stable DGA
isomorphism induced by a type II move is necessary.  Suppose that
$(\alg', \df')$ is the DGA for a knot diagram before a type II move
and that $(\alg, \df)$ is the DGA afterward.  As shown in \cite{chv},
the type II move gives rise to a DGA isomorphism
\begin{equation*}
  \psi: (\alg, \df) \to S(\alg', \df').
\end{equation*}
In particular, note that this map preserves grading.  If $a$ and $b$
are the two new generators that appear during a type II move, then the
first step in defining $\psi$ is to order the generators of \alg\ by
height: let $\{x_1, \ldots, x_N\}$ denote generators of height less
than that of $a$ in increasing height order and let $\{y_1, \ldots,
y_M\}$ denote generators of height greater than that of $b$ in
increasing height order.  Note that, since \df\ lowers height, $\df
y_j$ does not contain any generators $y_k$ with $k \geq j$.

It is possible to construct a dip in the plat diagram so that this
ordering takes on the following form.  Suppose the strand $k$ is
pushed over strand $l$.  Each $x_j$ either lies to the left of the
dip, or $x_j = a_{nm}$ or $b_{nm}$ with $n-m \leq k-l$.  Similarly,
$y_j$ either lies to the right of the dip, or $y_j = a_{nm}$ or
$b_{nm}$ with $n- m > k-l$.

The definition of the map $\psi$ needs a \emph{vector space} map
$H$ defined on $S(\alg')$ by:
\begin{equation*} 
  H(w) = \begin{cases}
    0 & w \in \alg' \\
    0 & w=Q\beta R \quad \text{with\ } Q \in \alg' \\
    Q\beta R & w=Q\alpha R \quad \text{with\ } Q \in \alg'.
  \end{cases}
\end{equation*}  

Also write $\df b = a + v$, where $v$ is a sum of words consisting
entirely of the letters $x_1, \ldots, x_N$.  Inductively define maps
$\psi_i$ on the generators of $\alg$ by:
\begin{equation*}
  \psi_0 (w) = \begin{cases}
    \beta & w = b \\
    \alpha + v & w = a \\
    w & \text{otherwise}
  \end{cases}
\end{equation*}
and
\begin{equation*}
  \psi_i (w) = \begin{cases}
    y_i + H\psi_{i-1}(\df y_i) & w=y_i \\
    \psi_{i-1}(w) & \text{otherwise.}
  \end{cases}
\end{equation*}
That the resulting map $\psi = \psi_M$ is a DGA isomorphism between
$\alg$ and $S(\alg')$ was proven in \cite{chv}.\footnote{This appears
  to be slightly different from the map given in \cite{chv, ens}; it
  is not hard to check, however, that the definition is equivalent.}

If there is an augmentation $\varepsilon'$ on $S(\alg')$, then
$\varepsilon = \varepsilon' \psi$ is an augmentation on \alg.  It is
straightforward to see that $\varepsilon(x_j) = \varepsilon'(x_j)$ and
that:
\begin{equation} \mylabel{eqn:trans-aug}
  \varepsilon(b) = \varepsilon'(\beta) \text{ and } \varepsilon(a) =
  \varepsilon'(v).
\end{equation}
Recall that if $|\beta| = 0$, then $\varepsilon'(\beta)$ may be chosen
arbitrarily.  In a plat diagram, there is a straightforward inductive
condition to determine if $\varepsilon$ will differ from
$\varepsilon'$ on a generator $y_j$:

\begin{lem} \mylabel{lem:y-condition}
  After a type II move involved in making a dip in a plat diagram,
  suppose that $\varepsilon(y_i)$ has been determined for all $i < j$.
  Then $\varepsilon'(y_j) \neq \varepsilon(y_j)$ if and only if
  $\varepsilon'(\beta) =1$ and there exists an odd number of terms in
  $\df y_j$ that are of the form $QaR$, where $Q, R \in \alg'$,
  $\varepsilon(Q) = 1$ and $\varepsilon(R) = 1$.
\end{lem}

\begin{proof}
  Since 
  \begin{equation} \mylabel{eqn:psi-yj}
    \psi(y_j) = y_j + H\psi(\df y_j),
  \end{equation}
  the augmentations $\varepsilon$ and $\varepsilon'$ disagree on $y_j$
  if and only if $\varepsilon'(H \psi (\df y_j)) \neq 0$. The proof
  that the latter is equivalent to the second condition in the lemma
  proceeds by induction on $j$.
  
  For $j=1$, let $P$ be the sum of terms in $\df y_1$ that do not
  contain $a$.  Lemma~\ref{lem:dip-disks} implies that
  $\df y_1$ has the form:
  \begin{equation} \mylabel{eqn:df-yj}
    \df y_1 = P + \sum_k Q_k a R_k,
  \end{equation}
  where $Q_k, R_k \in \alg'$.  Since the differential lowers height,
  $P$ lies in the algebra generated by $\{x_1, \ldots, x_N, b\}$.  It
  follows that:
  \begin{equation*}
    \begin{split}
      H \psi(\df y_1) &= H\bigl(\psi(P) + \sum_k Q_k (\alpha + v) R_k \bigr) \\
      &= \sum_k Q_k \beta R_k,
    \end{split}
  \end{equation*}
  since the $Q_k \alpha R_k$ are the only terms containing $\alpha$.
  The lemma follows in this case.  This argument also shows that
  $\alpha$ does not appear in $\psi(y_1)$.
  
  In general, write out $\df y_j$ as in equation (\ref{eqn:df-yj}).
  As before, the generator $a$ only appears where indicated, and $Q_k$
  and $R_k$ lie in the algebra generated by $\{x_1, \ldots, x_N, y_1,
  \ldots, y_{j-1} \}$.  Inductively, $\psi(y_i)$ does not contain
  $\alpha$ for $i < j$, so the images of $Q_k$, $R_k$, and $P$ under
  $\psi$ do not contain $\alpha$. This implies that $H \psi (Q_k
  \alpha) = \psi(Q_k) \beta$.  Computing as before, then,
  \begin{equation} \mylabel{eqn:H-psi}
    H \psi(\df y_j) = \sum_k \psi(Q_k) \beta \psi(R_k).
  \end{equation}
  Once again, this implies that $\alpha$ does not appear in
  $\psi(y_j)$, so this fact may be used inductively.  The lemma now
  follows from (\ref{eqn:psi-yj}), (\ref{eqn:H-psi}), and the fact
  that $\varepsilon' \psi = \varepsilon$.
\end{proof}

\subsection{Extension of the Ruling}
\mylabel{ssec:tracing}

The heart of the proof of Theorem~\ref{thm:main} extends ruling paths
that start at a common left cusp over successive crossings to the
right. In the Lagrangian projection that comes from resolving a plat
diagram, label the crossings that correspond to crossings of the plat
by $q_1, \ldots, q_n$.  The extension procedure has three parts:
First, extend the ruling over $q_j$; then place a dip between $q_j$
and $q_{j+1}$; and finally construct an augmentation
$\varepsilon_{j+1}$ on the DGA of the newly dipped diagram.  The
augmentations $\varepsilon_{j+1}$ will have the following property:

\begin{pty}[R]
  At any dip, $a_{kl}$ is augmented if and only if the strands $k$ and
  $l$ are paired in the portion of the ruling between $q_j$ and
  $q_{j+1}$.
\end{pty}

The construction begins at the left cusps, where any ruling must pair
paths incident to the same cusp.  The first step is to construct
$\varepsilon_1$ on the diagram that results from placing a dip between
the left cusps and $q_1$.  Consider the type II move that pushes
strand $k$ over strand $l$, and use the notation for augmentations and
generators that was set up around equation (\ref{eqn:trans-aug}).
There are three considerations that go into computing $\varepsilon$
from $\varepsilon'$:
\begin{enumerate}
\item A choice for $\varepsilon'(\beta)$ must be made.  In this case,
  choose $\varepsilon'(\beta) = 0$; it immediately follows from
  (\ref{eqn:trans-aug}) that $\varepsilon(b_{kl}) = 0$.
\item The value of $\varepsilon(a_{kl})$ is determined from
  $\varepsilon'(v_{kl})$ via (\ref{eqn:trans-aug}).  In this case,
  Figure~\ref{fig:left-cusp-dip} shows that $v_{kl}$ is a sum of words
  in $b_{ij}$ (for $(i,j) \prec (k,l)$) and contains a $1$ if $(k,l) =
  (2m, 2m+1)$ for some $m$.  Since $\varepsilon'(b_{ij}) = 0$ for all
  $(i,j) \prec (k,l)$ by step (1), it is simple to compute
  $\varepsilon'(v_{kl})$, and hence $\varepsilon(a_{kl})$:
  \begin{equation} \label{eqn:a-kl-augm}
    \varepsilon(a_{kl}) = \varepsilon'(v_{kl}) = \begin{cases}
      1 & (k,l) = (2m, 2m+1) \\
      0 & \text{otherwise.}
    \end{cases}
  \end{equation}

  \begin{figure}
    \centerline{\includegraphics{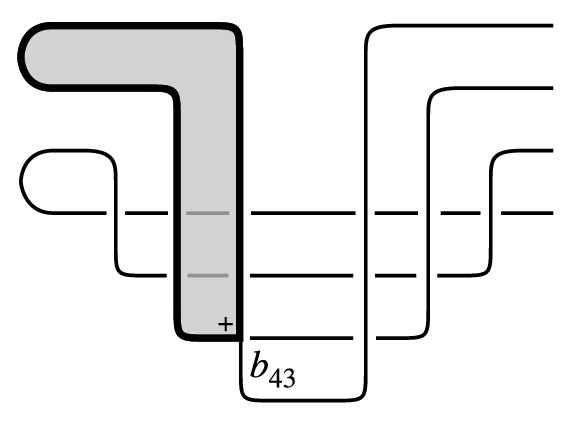}}
    \caption{The dip next to the left cusp with an augmented disk in $v_{43}$.}  
    \mylabel{fig:left-cusp-dip}
  \end{figure}
  
\item Finally, Lemma~\ref{lem:y-condition} is used to check if there
  are any ``corrections'' to other $a_{ij}$ generators with $(i,j)
  \prec (k,l)$ but $i-j \geq k-l$.  In this case, since
  $\varepsilon'(\beta) = 0$, no such changes can occur.
\end{enumerate}
At all stages, then, this process gives an augmentation that satisfies
(\ref{eqn:a-kl-augm}), and hence $\varepsilon_1$ satisfies Property
(R).

Now begin the extension procedure proper.  At the crossing $q_j$,
extend the ruling paths as follows: if $\varepsilon_j(q_j) = 1$ and
the ruling to the left of $q_j$ matches the situation in
configurations (a), (b), or (c) in Figure~\ref{fig:switch-config},
then there is a switch at $q_j$.  Otherwise, there is no switch. By
construction, the ruling paths have only (graded or $\rho$-graded)
normal switches.

The next part of the extension procedure is to understand the
augmentation $\varepsilon_{j+1}$ that results from the construction of
a dip between $q_j$ and $q_{j+1}$ using the three steps above.  The
choice of augmentations on the $\beta$ generators in step (1) should
lead to $\varepsilon_{j+1}$ satisfying Property (R) if $\varepsilon_j$
does.  The exact choice of augmentations depends on
$\varepsilon_j(q_j)$ and the configuration of the ruling near the
crossing $q_j$.

First, suppose that $\varepsilon_j(q_j) = 0$ and consider the Type II
move that pushes strand $k$ over strand $l$.  For step (1), choose
$\varepsilon'(\beta) = 0$.
  
For step (2), consider $\varepsilon'(v_{kl})$.  Since neither $q_j$,
nor any crossing in the $b$-lattice, is augmented, the only totally
augmented disks in $v_{kl}$ have a positive corner at $b_{kl}$ and a
single augmented negative corner in the $a$-lattice to the left of
$q_j$; see Figure~\ref{fig:eps0-disks}. If such a disk exists, the
negative corner must occur where two ruling strands cross each other,
since $\varepsilon'$ satisfies property (R) on the $a$-lattice to the
left.  The facts that $q_j$ is not switched in the ruling and that
there are no other corners on the disk imply that $b_{kl}$ --- and
hence $a_{kl}$ --- must also be crossings of ruling strands.  Thus,
$\varepsilon'(v_{kl}) = \varepsilon(a_{kl}) = 1$ if and only if $k$
and $l$ are paired in the ruling.
  
Finally, since $\varepsilon'(\beta) = 0$, Lemma~\ref{lem:y-condition}
shows that there are no corrections to the augmentations of $a_{ij}$
for $(i,j) \prec (k,l)$.  Thus, the previous paragraph shows that
$\varepsilon_{j+1}$ satisfies property (R).

\begin{figure}
  \centerline{\includegraphics{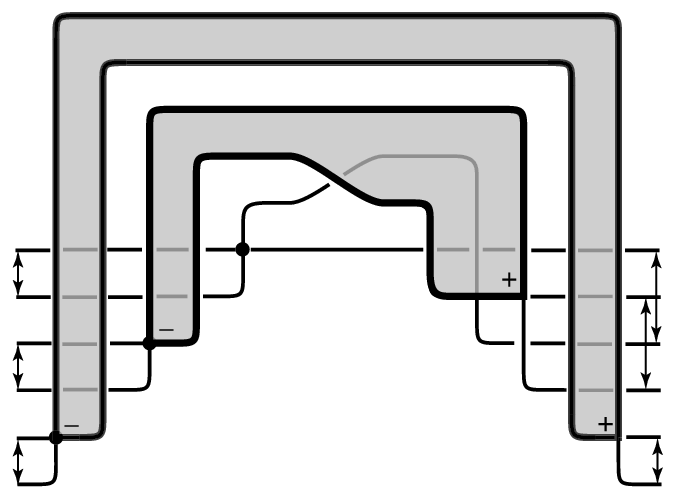}}
  \caption{In the case that $\varepsilon_j(q_j) = 0$, the only totally augmented disks in $v_{kl}$ have a positive corner at
    $b_{kl}$ and a single augmented negative corner in the $a$-lattice
    to the left of $q_j$.  The filled dots represent augmented
    corners, and the arrows represent strands paired in the ruling.}
  \mylabel{fig:eps0-disks}
\end{figure}
  
From now on, assume that $\varepsilon_j(q_j) = 1$; the proof will
examine each configuration in Figure~\ref{fig:switch-config} in turn.
For configuration (a), suppose that the strands $i$ and $i+1$ cross
and that these strands are paired with $L$ and $K$, respectively.
That is, $K > i+1 >i>L$.  Divide the dipping process into three parts:

\begin{description}
\item[$(k,l) \prec (i+1,i)$] Choose $\varepsilon'(\beta)=0$.  To
  determine $\varepsilon(a_{kl})$, consider totally augmented disks in
  $v_{kl}$.  As before, the leftmost negative corner of a totally
  augmented disk must involve strands paired in the ruling. If neither
  $k$ nor $l$ is a crossing strand, then, as above,
  $\varepsilon'(a_{kl}) = \varepsilon(a_{kl}) = 1$ if and only if $k$
  and $l$ are paired in the ruling.  Otherwise, Figure~\ref{fig:config-a-1}
  shows that there is one totally augmented disk in each of
  $v_{i+1,L}$ and $v_{i,L}$.  Thus,
  \begin{equation} \mylabel{eqn:config-a-1}
    \varepsilon(a_{i+1,L}) = \varepsilon(a_{iL}) = 1.
  \end{equation}
  
  \begin{figure}
    \centerline{\includegraphics{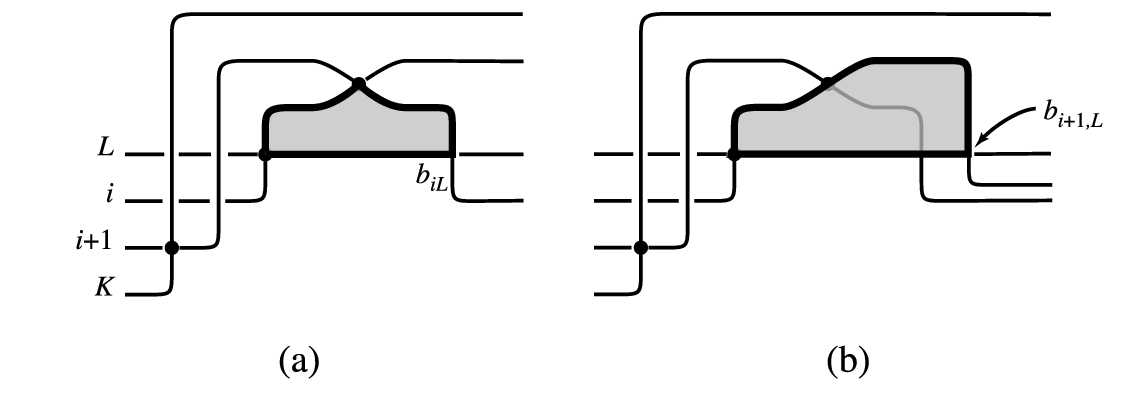}} 
    \caption{The
      totally augmented disks in (a) $v_{iL}$ and (b) $v_{i+1,L}$ in the
      proof of property (R) for configuration (a).}
    \mylabel{fig:config-a-1}
  \end{figure}
  
  Since $\varepsilon'(\beta)=0$, there are no corrections to the
  augmentations of previously constructed crossings in the
  $a$-lattice.
  
\item[$(k,l) = (i+1,i)$] First, note that $|b_{i+1,i}|=0$ if the
  augmentation is graded: the Maslov indices of the crossing strands
  must agree, and $b_{i+1,i}$ involves the crossing strands.  A
  similar fact holds for a $\rho$-graded augmentation.  Hence, it is
  possible to choose $\varepsilon'(\beta)=1$; it follows that
  $\varepsilon(b_{i+1,i}) = 1$ as well.
  
  It is easy to see that $v_{i+1,i} = 0$, so $\varepsilon(a_{i+1,i}) =
  0$.  There is one correction to consider.  The disk in
  Figure~\ref{fig:config-a-2} contributes the term $a_{i+1,i} a_{iL}$
  to $\df a_{i+1,L}$.  This is the only disk with a negative corner at
  $a_{i+1,i}$ whose other negative corners are augmented since
  $a_{iL}$ is the only crossing involving strand $L$ that is
  augmented.  Equation (\ref{eqn:config-a-1}) shows that
  $\varepsilon(a_{i+1,L})=\varepsilon(a_{iL})=1$, so
  Lemma~\ref{lem:y-condition} implies that
  \begin{equation} \mylabel{eqn:config-a-2}
      \varepsilon(a_{i+1,L}) = 0.
    \end{equation}
    Thus, the augmentation on all crossings created up to this point
    satisfies property (R).

  \begin{figure}
    \centerline{\includegraphics{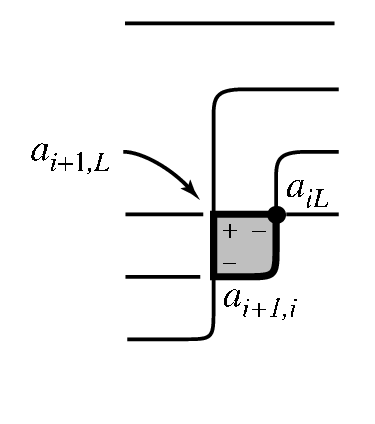}} 
    \caption{The correction disk 
      in the proof of property (R) for configuration (a).}
    \mylabel{fig:config-a-2}
  \end{figure}
  
\item[$(k,l) \succ (i+1,i)$] Choose $\varepsilon'(\beta) = 0$. As in
  the case of $(k,l) \prec (i+1,i)$, if neither strand is a crossing
  strand, then the augmentation for $a_{kl}$ matches the augmentation
  in the $a$-lattice to the left. On the other hand,
  Figure~\ref{fig:config-a-3} shows that there is a single totally
  augmented disk in $v_{K,i+1}$ and two totally augmented disks in
  $v_{Ki}$.  Thus,
    \begin{equation} \mylabel{eqn:config-a-3}
      \varepsilon(a_{K,i+1}) = 1 \quad \text{and} \quad \varepsilon(a_{Ki}) = 0.
    \end{equation}

  \begin{figure}
    \centerline{\includegraphics{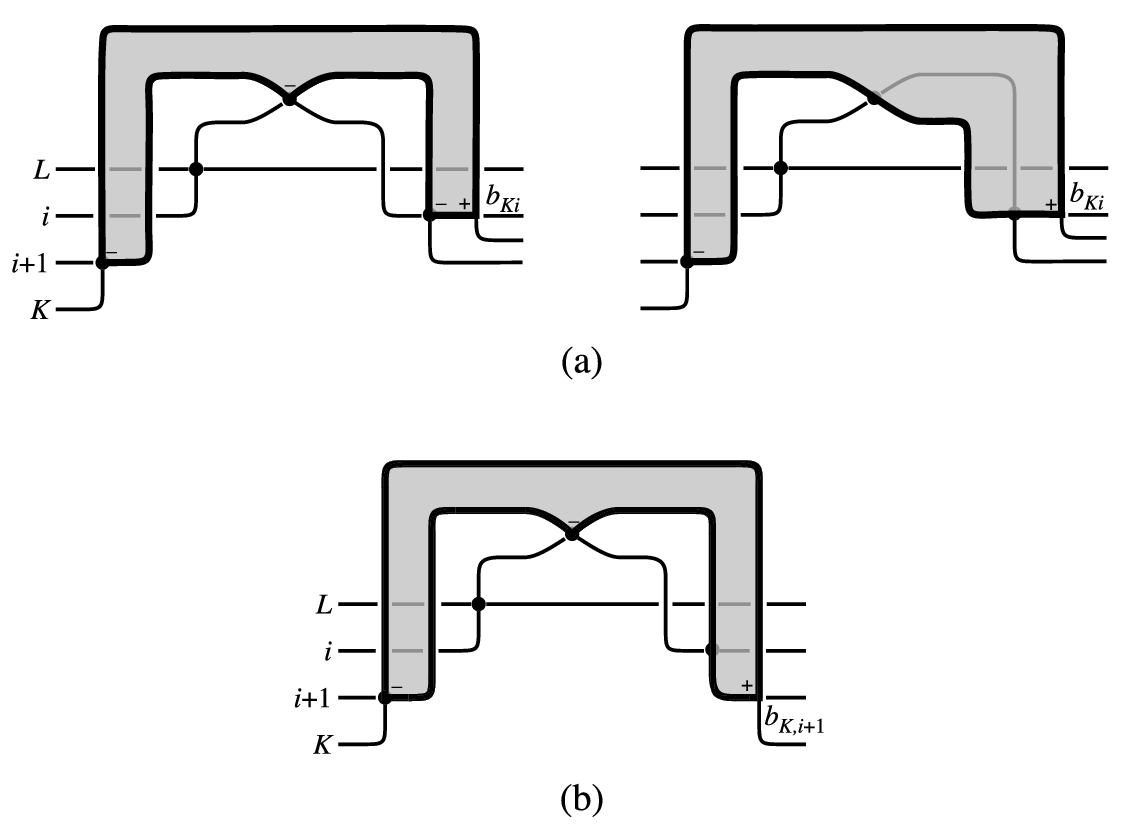}} 
    \caption{(a) The two totally augmented disks in $v_{Ki}$ and (b) the totally augmented disk in $v_{K,i+1}$
      in the proof of property (R) for configuration (a).}
    \mylabel{fig:config-a-3}
  \end{figure}
  
  Since $\varepsilon'(\beta)=0$, there are no corrections to the
  augmentations of previously constructed crossings in the
  $a$-lattice.
\end{description}

The end result is an augmentation that satisfies property (R) on the
new $a$-lattice: for crossing strands, equations
(\ref{eqn:config-a-1}, \ref{eqn:config-a-2}, \ref{eqn:config-a-3})
show that only $a_{iL}$ and $a_{K,i+1}$ are augmented; otherwise, the
augmentation is simply transferred from the $a$-lattice to the left.


  
The next case to consider is configuration (b).  Again, suppose that
the crossing strands are $i$ and $i+1$, paired with $K$ and $L$,
respectively, so that $i+1>i>K>L$.  This time, the dipping process
should be divided into five steps:

\begin{description}
\item[$(k,l) \prec (K,L)$] As in the first case in configuration (a),
  set $\varepsilon'(\beta)=0$ and transfer the augmentations from the
  $a$-lattice on the left.
  
\item[$(k,l) = (K,L)$] Note that $|b_{KL}| = 0$ for a graded
  augmentation: the crossing strands have the same Maslov index, and
  hence so do the companion strands since they both lie below their
  corresponding crossing strands.  Thus, it is possible to set
  $\varepsilon'(\beta) = 1$, and hence obtain $\varepsilon(b_{KL})=1$.
  
  Since $K$ and $L$ are not paired in the ruling and are not crossing
  strands, $\varepsilon'(v_{KL}) = 0$, so $\varepsilon(a_{KL}) = 0$.
  Further, there are no corrections, as any disk in the $a$-lattice
  with a negative corner at $a_{KL}$ must have an augmented negative
  corner of the form $a_{L*}$ (see Figure~\ref{fig:config-a-2}).
  Since $L$ is paired with $i$, the only augmented crossing of this
  form has yet to appear in the dip.

\item[$(K,L) \prec (k,l) \prec (i+1,i)$] Set $\varepsilon'(\beta) =
  0$.  There are several augmented disks contributing to $v_{kl}$; see
  Figure~\ref{fig:config-b-1}:
  \begin{itemize}
  \item Two for $v_{iL}$, and hence $\varepsilon(a_{iL}) = 0$. Note
    that one of these disks uses the fact that $\varepsilon(b_{KL}) =
    1$.
  \item One for $v_{iK}$, and hence $\varepsilon(a_{iK}) = 1$.
  \item One for $v_{i+1,L}$, and hence $\varepsilon(a_{i+1,L}) = 1$.
    Note that the existence of this disk relies on the fact that
    $\varepsilon(b_{KL}) = 1$.
  \item One for $v_{i+1,K}$, and hence $\varepsilon(a_{i+1,K}) = 1$.
  \end{itemize}
  Since $\varepsilon'(\beta) = 0$, there are no corrections at this stage.
  
  \begin{figure}
    \centerline{\includegraphics{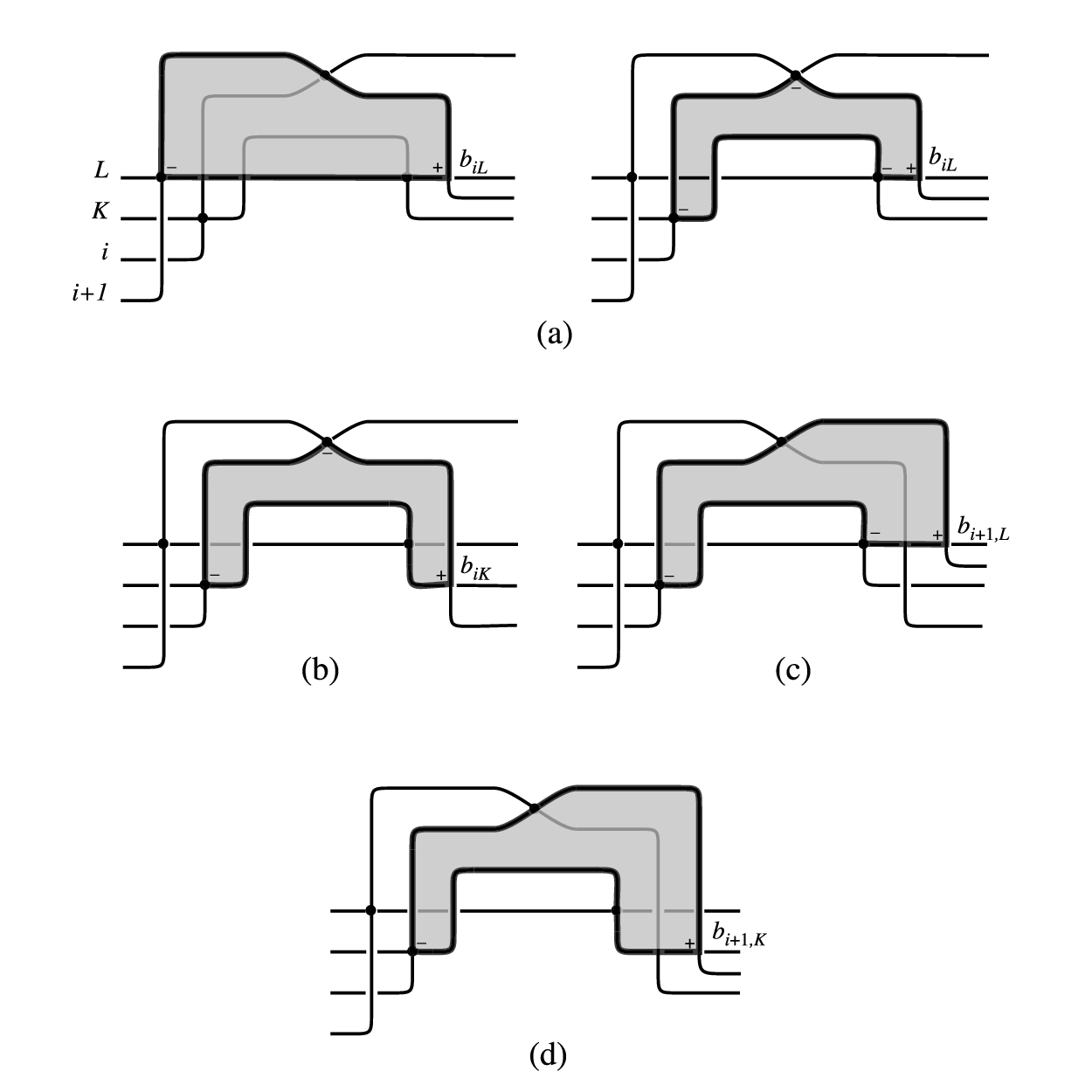}} 
    \caption{The totally augmented disks in (a) $v_{iL}$, (b)
      $v_{iK}$, (c) $v_{i+1,L}$, and (d) $v_{i+1,K}$ in the proof of
      property (R) for configuration (b).}  \mylabel{fig:config-b-1}
  \end{figure}
  
\item[$(k,l) = (i+1,i)$] Set $\varepsilon'(\beta) = 1$.  As usual,
  $\varepsilon'(v_{i+1,i}) = 0$, so $\varepsilon(a_{i+1,i}) = 0$.
  There is one correction in this case: since one term in $\df
  a_{i+1,K}$ is $a_{i+1,i} a_{iK}$, Lemma~\ref{lem:y-condition}
  implies that $\varepsilon(a_{i+1,K})$ changes to $0$.


\item[$(k,l) \succ (i+1,i)$] Set $\varepsilon'(\beta)=0$.  As in the
  final case in configuration (a), the augmentation is simply
  transferred from the the dip on the left.
\end{description}

In sum, the augmentation on the new dip satisfies property (R): for
crossing strands, only $a_{i+1,L}$ and $a_{iK}$ are augmented;
otherwise, the augmentation is simply transferred from the $a$-lattice
to the left.

The arguments for the other configurations is similar; see
Table~\ref{tbl:config-augm} for a list of which $\beta$ generators to
augment in each case.  This completes the extension of the ruling and
of $\varepsilon_j$ over a dip.

\begin{table}

  \begin{tabular}{c|c}
    Configuration & Augmented Generators \\ \hline \hline
    a & Crossing \\ \hline
    b,c & Crossing and companion \\
    \hline
    d & None \\ \hline
    e,f & Companion
  \end{tabular}

  \caption{Which $\beta$ generators are augmented in each
    configuration (with $\varepsilon_j(q_j) = 1$)?}
  \mylabel{tbl:config-augm}
\end{table}

As mentioned above, the proof of Theorem~\ref{thm:main} will be
complete if the paired ruling paths match at the right cusps.  This
holds true if and only if, in the dip just to the left of the right
cusps, $a_{2k, 2k-1}$ is augmented for $k=1, \ldots, c$.  As shown in
Figure~\ref{fig:right-cusp-dip}, the differential of the $k^{th}$
right cusp in the dipped diagram is:
\begin{equation*}
  \df q_{n+k} = 1 + a_{2k, 2k-1}.
\end{equation*}
Since $\varepsilon_1$ satisfies property (R), the inductive extension
argument above shows that $\varepsilon_n$ does as well. The fact that
$\varepsilon_n$ is a genuine augmentation implies that $a_{2k, 2k-1}$
is augmented.  Theorem~\ref{thm:main} follows since $\varepsilon_n$
obeys property (R).

\begin{figure}
  \centerline{\includegraphics{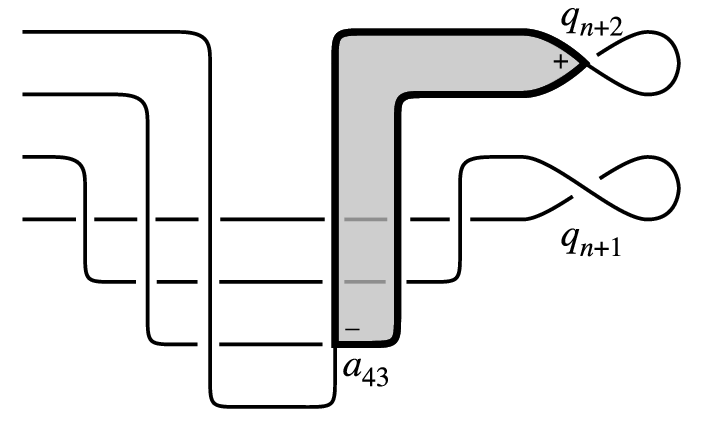}}
  \caption{The dip next to the right cusp with the disk that gives the
    $a_{21}$ term in $\df q_{n+1}$.}  
  \mylabel{fig:right-cusp-dip}
\end{figure}

\begin{rem}
  The proof can be refined to give an algorithm for constructing a
  ruling from the augmentation, and can even be carried out without
  passing to the dipped diagram.  As in the proof, the idea is to
  extend the ruling over a crossing $q_j$ given the value
  $\varepsilon_j(q_j)$.  Before, it was not necessary to explicitly
  find these values, but it \emph{is} possible to determine them.
  
  The key to finding $\varepsilon_j(q_k)$ for $k > j$ is
  Lemma~\ref{lem:y-condition}.  Disks of the form $QaR$, where $Q, R
  \in \alg'$, appear in the original plat as disks with a positive
  corner at $q_k$, negative corners at $Q$ and $R$, and a line segment
  to the right of $q_j$ that joins the crossing strands (if the
  $\beta$ generator at the crossing is augmented) or the companion
  strands (if the corresponding $\beta$ generator is augmented); see
  Figure~\ref{fig:crossing-ends}.  The value of $\varepsilon_j(q_k)$
  differs from that of $\varepsilon_{j+1}(q_k)$ if there is an odd
  number of these disks with $\varepsilon_{j+1}(Q)=1$ and
  $\varepsilon_{j+1}(R) = 1$.  If two $\beta$ generators are
  augmented, then the procedure should be performed once for each
  $\beta$ with the $\beta$ between the lower-numbered strands going
  first.
  
  Figure~\ref{fig:ruling-proc} demonstrates the procedure on the
  trefoil with the augmentation that marks all three crossings of
  degree 0.  Note that all of the disks involved in adjusting the
  augmentation have no negative corners, so $Q$ and $R$ are always $1$
  in Lemma~\ref{lem:y-condition}, and the condition is easy to apply.
\end{rem}

\begin{figure}
  \centerline{\includegraphics{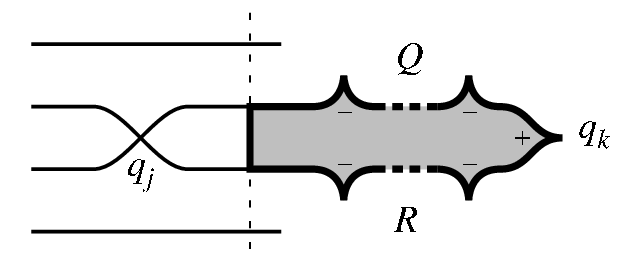}}
  \caption{A schematic picture of a disk that contributes to changing 
    $\varepsilon_{j+1}(q_k)$ in configuration (a).}
  \mylabel{fig:crossing-ends}
\end{figure}

\begin{figure}
  \centerline{\includegraphics{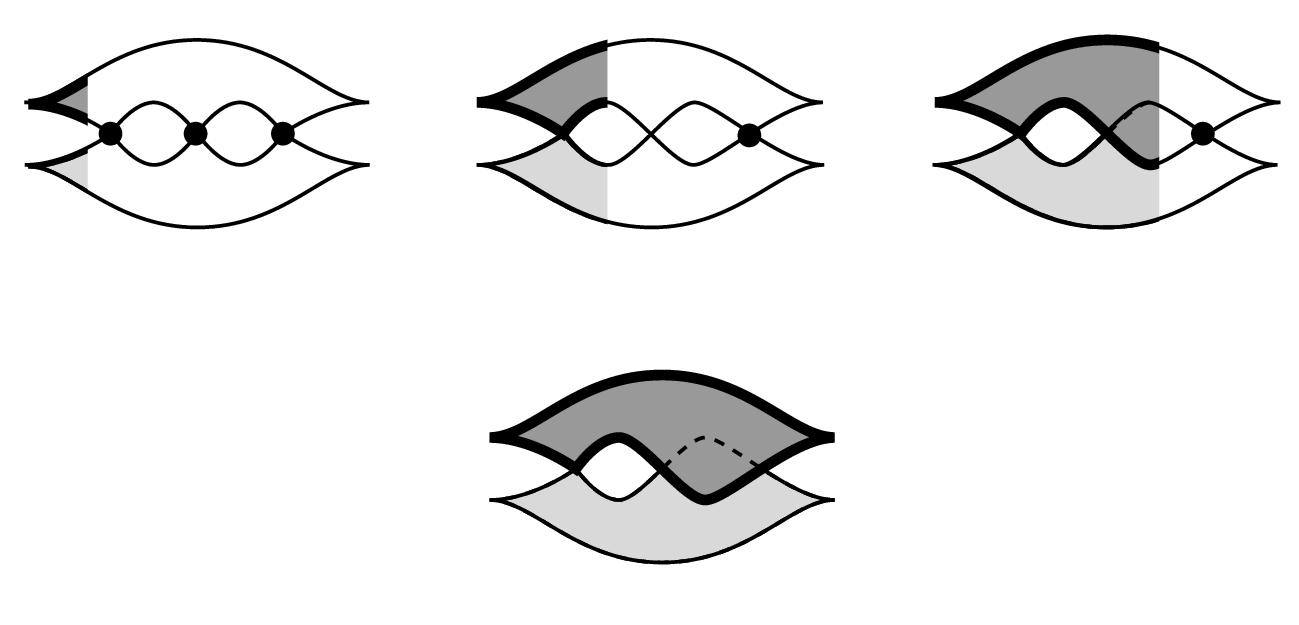}}
  \caption{The result of carrying out the procedure in this section on 
    the trefoil with the augmentation $\varepsilon_3$ that sends all
    three crossings of degree 0 to $1$.}  \mylabel{fig:ruling-proc}
\end{figure}

\section{Rotation Number and Rulings}
\mylabel{sec:augm-r}

This brief section contains the proof of Theorem~\ref{thm:rot-ruling}.
By Theorem~\ref{thm:main}, it suffices to prove that if an oriented
front diagram of $K$ has a $2$-graded normal ruling then $r(K)=0$.
  
It is easy to check that the strands at a crossing with even grading
are both oriented to the left or both to the right.  This implies that
the boundary of a disk in a graded normal ruling inherits a coherent
orientation from the knot, and hence that each disk pairs an upward
(resp.  downward) right cusp with a downward (resp. upward) left cusp.
Thus,
\begin{eqnarray*}
  2 r(K) &=& \#\text{down cusps} - \#\text{up cusps} \\
  &=& \#\text{down right cusps} - \#\text{up left cusps} \\ 
  & & \quad +
  \#\text{down left cusps} - \#\text{up right cusps} \\
  &=& 0.
\end{eqnarray*}

\bibliographystyle{amsplain}
\bibliography{rulings}

\providecommand{\bysame}{\leavevmode\hbox to3em{\hrulefill}\thinspace}
\providecommand{\MR}{\relax\ifhmode\unskip\space\fi MR }
\providecommand{\MRhref}[2]{%
  \href{http://www.ams.org/mathscinet-getitem?mr=#1}{#2}
}
\providecommand{\href}[2]{#2}
\begin{thebibliography}{10}

\bibitem{bennequin}
D.~Bennequin, \emph{Entrelacements et equations de {P}faff}, Asterisque
  \textbf{107--108} (1983), 87--161.

\bibitem{chv}
Yu. Chekanov, \emph{Differential algebra of {L}egendrian links}, Invent. Math.
  \textbf{150} (2002), 441--483.

\bibitem{chv:survey}
\bysame, \emph{Invariants of {L}egendrian knots}, Proceedings of the
  International Congress of Mathematicians, Vol. II (Beijing, 2002) (Beijing),
  Higher Ed. Press, 2002, pp.~385--394.

\bibitem{chv-pushkar}
Yu. Chekanov and P.~Pushkar, \emph{The combinatorics of fronts of {L}egendrian
  knots}, Preprint., 2004.

\bibitem{ees:graph-trees}
T.~Ekholm, J.~Etnyre, and M.~Sullivan, In preparation.

\bibitem{yasha-fraser}
Y.~Eliashberg and M.~Fraser, \emph{Classification of topologically trivial
  {L}egendrian knots}, Geometry, topology, and dynamics (Montreal, PQ, 1995),
  Amer. Math. Soc., Providence, RI, 1998, pp.~17--51.

\bibitem{egh}
Y.~Eliashberg, A.~Givental, and H.~Hofer, \emph{Introduction to symplectic
  field theory}, Geom. Funct. Anal. (2000), no.~Special Volume, Part II,
  560--673, GAFA 2000 (Tel Aviv, 1999).

\bibitem{etnyre:knot-intro}
J.~Etnyre, \emph{{L}egendrian and transversal knots}, To appear in the Handbook
  of Knot Theory, 2003.

\bibitem{etnyre-honda:knots}
J.~Etnyre and K.~Honda, \emph{Knots and contact geometry}, J. Symplectic Geom.
  \textbf{1} (2002), no.~1, 63--120.

\bibitem{ens}
J.~Etnyre, L.~Ng, and J.~Sabloff, \emph{Invariants of {L}egendrian knots and
  coherent orientations}, J. Symplectic Geom. \textbf{1} (2002), no.~2,
  321--367.

\bibitem{fuchs:augmentations}
D.~Fuchs, \emph{Chekanov-{E}liashberg invariant of {L}egendrian knots:
  existence of augmentations}, J. Geom. Phys. \textbf{47} (2003), no.~1,
  43--65.

\bibitem{fuchs-ishk}
D.~Fuchs and T.~Ishkhanov, \emph{Invariants of {L}egendrian knots and
  decompositions of front diagrams}, Moscow Math. J. (2004), To appear.

\bibitem{fukaya-oh}
K.~Fukaya and Y.-G. Oh, \emph{Zero-loop open strings in the cotangent bundle
  and {M}orse homotopy}, Asian J. Math. \textbf{1} (1997), no.~1, 96--180.

\bibitem{lenny:computable}
L.~Ng, \emph{Computable {L}egendrian invariants}, Topology \textbf{42} (2003),
  no.~1, 55--82.

\bibitem{lenny:knot-invts-1}
\bysame, \emph{Knot and braid invariants from contact homology i}, Available on
  arXiv as math.GT/0302099, 2003.

\bibitem{lenny-lisa}
L.~Ng and L.~Traynor, \emph{Legendrian solid-torus links}, Available on arXiv
  as math.SG/0407068, 2004.

\bibitem{lecnotes}
J.~Sabloff, \emph{Invariants for {L}egendrian knots from contact homology}, In
  preparation.

\bibitem{lisa:links}
L.~Traynor, \emph{Generating function homology for {L}egendrian links}, Geom.
  and Top. \textbf{5} (2001), 719--760.

\bibitem{zhu}
Ke~Zhu, \emph{Degeneration of the moduli space of j-holomorphic discs and
  legendrian contact homology}, Ph.D. thesis, Stanford University, 2004.

\end{thebibliography}

\end{document}